\documentclass[dvipdfmx]{amsart}
\usepackage{amsmath}
\usepackage{amssymb} 
\usepackage{amscd} 
\usepackage{graphicx} 
\usepackage{epsfig}
\usepackage[dvips]{color}

\newtheorem{theorem}{Theorem}[section]
\newtheorem{lemma}[theorem]{Lemma}

\newtheorem{remark}[theorem]{Remark}
\newtheorem{question}[theorem]{Question}

\theoremstyle{definition}

\numberwithin{equation}{section}
\numberwithin{figure}{section}
\numberwithin{table}{section}

\begin{document}
\baselineskip 14pt

\title{Hyperbolic, L--space knots and exceptional Dehn surgeries}

\author[K. Motegi]{Kimihiko Motegi}
\address{Department of Mathematics, Nihon University, 
3-25-40 Sakurajosui, Setagaya-ku, 
Tokyo 156--8550, Japan}
\email{motegi@math.chs.nihon-u.ac.jp}
\thanks{
The author has been partially supported by JSPS Grants--in--Aid for Scientific 
Research (C), 26400099, The Ministry of Education, Culture, Sports, Science and Technology, Japan and Joint Research Grant of Institute of Natural Sciences at Nihon University for 2014. 
}

\author[K. Tohki]{Kazushige Tohki}
\address{Graduate School of Integrated Basic Sciences, Nihon University, 
3-25-40 Sakurajosui, Setagaya-ku, 
Tokyo 156--8550, Japan}

\dedicatory{}

\begin{abstract}
A knot in the $3$--sphere is called an \textit{L--space knot} if it admits a nontrivial Dehn surgery 
yielding an L--space.  
Like torus knots and Berge knots, 
many L--space knots admit also a Seifert fibered surgery. 
We give a concrete example of a hyperbolic, L--space knot which has no exceptional surgeries, 
in particular, no Seifert fibered surgeries.  
\end{abstract}

\maketitle

{
\renewcommand{\thefootnote}{}
\footnotetext{
}

\section{Introduction}
\label{section:Introduction}

For any rational homology $3$--sphere $M$, 
the rank of the Heegaard Floer homology $\widehat{\mathrm{HF}}(M)$ is bounded below by the order of 
$H_1(M; \mathbb{Z})$. 
If the rank of $\widehat{\mathrm{HF}}(M)$ is equal to the order of $H_1(M; \mathbb{Z})$, 
then $M$ is called an \textit{L--space}. 
The class of L--spaces includes lens spaces (except $S^2 \times S^1$), 
and more generally, $3$--manifolds with elliptic geometry \cite[Proposition~2.3]{OS3}. 
An efficient way to find L--spaces is through ``bootstrapping" a known ``L--space surgeries" on a knot. 
A Dehn surgery is called an \textit{L--space surgery} if the resulting $3$--manifold is an L--space, 
and a knot admitting nontrivial L--space surgery is called an \textit{L--space knot}. 
Since torus knots and Berge knots \cite{Berge2} admit surgeries yielding lens spaces, 
these are L--space knots. 
Ozsv\'ath-Szab\'o \cite[Proposition 9.6]{OS4}  (\cite[Lemma 2.13]{Hedden}) gives us 
a complete information about the set of L--space surgeries on an L--space knot: 

\vskip 0.3cm

$\bullet$ If $K$ is a nontrivial, L--space knot of genus $g(K)$, 
then $r$--surgery on $K$ results in an L--space if $r \ge 2g(K) - 1$ or $r \le -2g(K)+ 1$. 

\vskip 0.3cm

This result, together with Thurston's hyperbolic Dehn surgery theorem
\cite{T1, T2, BePe, PetPorti, BoileauPorti}, 
says that each hyperbolic, L--space knot produces infinitely many hyperbolic L--spaces by Dehn surgery. 
For instance, 
a hyperbolic Berge knot produces infinitely many hyperbolic L--spaces. 

\smallskip

Like torus knots and Berge knots, 
many L--space knots admit also a Seifert fibered surgery, 
i.e. a surgery yielding a Seifert fiber space. 
Among Montesinos knots, 
recent results of Lidman-Moore \cite{LM} and Baker-Moore \cite{BM} show that the only L--space knots are the pretzel knots 
$P(-2, 3, 2n+1)$ and the torus knots $T_{2n+1, 2}$ $(n \ge 0)$ and their mirror images, 
each of which admits a Seifert fibered surgery \cite{OS3, LM}. 
Furthermore, in \cite{Mote} one can find a large number of twist families of hyperbolic, 
L--space knots each of which admits also a Seifert fibered surgery. 
To the best of our knowledge, 
there is no explicitly known examples of hyperbolic, L--space knots which have no Seifert fibered surgeries, 
though we expect there should be many. 

\begin{question}
\label{Lspace=Seifert}
Does any hyperbolic, L--space knot admit also a Seifert fibered L--space surgery?
\end{question}

The aim of this note is to demonstrate:  

\begin{theorem}
\label{main}
There exists a hyperbolic, L--space knot which has no exceptional surgeries, 
in particular, no Seifert fibered surgeries.  
\end{theorem}

\begin{remark}
\label{satellite}
Using a cabling construction \cite[Theorem~1.10]{Hedden}, or more generally a satellite operation \cite[Theorem~1.3]{HLV}, 
we can obtain a satellite $($i.e. non-hyperbolic$)$ L--space knot which has no Seifert fibered surgeries. 
\end{remark}

\bigskip

\noindent
\textbf{Acknowledgments.}
We would like to thank Ken Baker for useful discussion, 
in which we learned effective utilization of almost alternating diagrams of the unknot in a study of L--space knots suggested by Josh Greene.  
We gratefully acknowledge Mario Eudave-Mu\~noz, Katura Miyazaki and Tatsuya Tsukamoto for useful conversations. 
Finally we would like to thank Neil Hoffman, Kazuhiro Ichihara and Hidetoshi Masai for their help with \texttt{fef.py} and HIKMOT. 

\bigskip

\section{Covering knots and Montesinos trick}

A \textit{tangle} $(B, t)$ is a pair of a $3$--ball $B$ and two disjoint arcs $t$ properly embedded in $B$. 
A tangle $(B, t)$ is trivial if 
there is a pairwise
homeomorphism from $(B, t)$ to $(D^2 \times I, \{x_1, x_2\}\times I)$, where $x_1, x_2$ are distinct points.
For tangles $(B, t)$ and $(B, t')$ with $\partial t = \partial t'$,
we say that they are \textit{equivalent} 
if there is a pairwise homeomorphism 
$h : (B, t) \to (B, t')$
satisfying $h|_{\partial B} =$ id. 

Let $U$ be the unit $3$--ball in $\mathbb{R}^3$,
and take $4$ points NW, NE, SE, SW on the boundary of $U$
so that
$\mathrm{NW} = (0, -\alpha, \alpha)$, 
$\mathrm{NE} = (0, \alpha, \alpha)$, 
$\mathrm{SE} = (0, \alpha, -\alpha)$, 
$\mathrm{SW} = (0, -\alpha, -\alpha)$, 
where $\alpha =\frac{1}{\sqrt{2}}$.
A tangle $(U, t)$ is a \textit{rational tangle}
if it is a trivial tangle with $\partial t =
\{\mathrm{NW, NE, SE, SW}\}$.
We can construct rational tangles from sequences of integers 
$a_1, a_2, \dots, a_n$ as shown in Figure~\ref{rtangle}, 
where the last horizontal twist $a_n$ may be $0$.
We consider that the tangle diagrams in Figure~\ref{rtangle}
is drawn on the $yz$--plane. 
Denote by $R(a_1, a_2, \dots, a_n)$ the associated rational tangle.  

\begin{figure}[htbp]
\begin{center}
\includegraphics[width=0.9\linewidth]{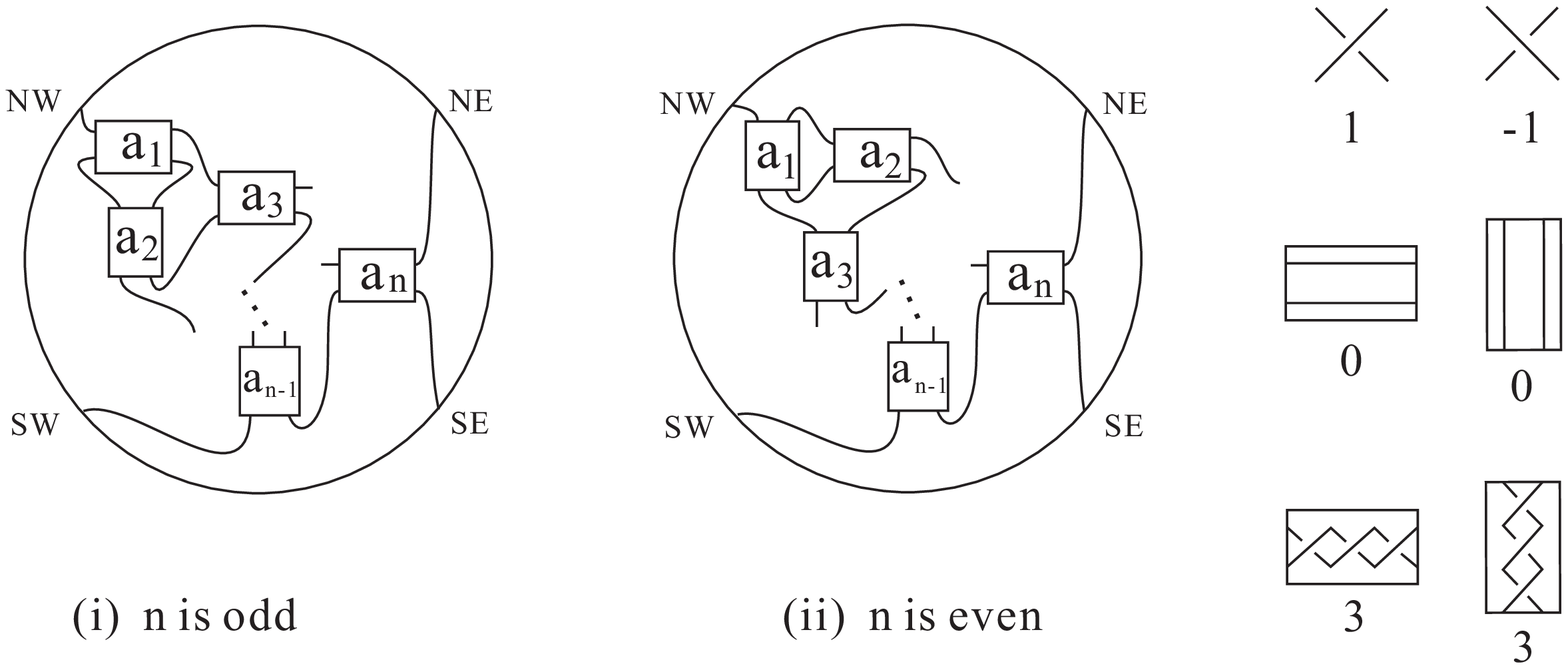}
\caption{Rational tangles}
\label{rtangle}
\end{center}
\end{figure}

Each rational tangle can be parametrized by 
$r \in \mathbb{Q} \cup \{ \infty\}$, 
where the rational number $r$ is given by
the continued fraction below. 
Thus we denote the rational tangle corresponding to $r$ by $R(r)$. 

$$ r\ =\ a_{n} + \cfrac 1 {a_{n-1}+ \cfrac 1 {
               \begin{array}{clr}
               \ & & \\[-5pt]
               \hspace{-25pt} \ddots & & \\[-10pt]
                      & \ \ \hspace{-5pt} a_2+ \cfrac{1}{a_1}
               \end{array}
               }}$$
               
Let $(U, t)$ be the rational tangle $R(\infty)$.
Considering $t$ is embedded in the $yz$--plane,
take the disk $D$ in the $yz$--plane such that
$\partial D$ is the union of $t$ and 
two arcs in $\partial U$: one connects
NW and NE, and the other connects SW and SE.
We call an arc in $D$ connecting the components of the interior of $t$
a \textit{spanning arc},
and the arc $D\cap \partial U$ connecting
NW and NE \textit{the latitude of $R(\infty)$}. 
See Figure~\ref{spanningarc}. 
The two--fold cover $\widetilde{U}$
of $U$ branched along $t$ is a solid torus.
Note that the preimages of the spanning arc
and the latitude are a core and a longitude $\lambda$
of the solid torus, respectively.
A meridian of a rational tangle $R(r)=(U, t')$ is
a simple closed curve in $\partial U -t'$ which bounds a disk
in $U -t'$ and a disk in $\partial U$ meeting $t'$
in two points.
Let $\mu_r (\subset \partial \widetilde{U})$
be a lift of a meridian of $R(r)$;
then $\mu_{\infty}$ is a meridian of the solid torus $\widetilde{U}$.
Furthermore, we note the following well-known fact.

\begin{lemma}
\label{-p/q}
Under adequate orientations we have
$[\mu_r] = -p[\mu_{\infty}] +q[\lambda]
\in H_1(\partial \widetilde{U})$,
where $r = \frac{p}{q}$ and $[\mu_{\infty}]\cdot[\lambda] =1$.
\end{lemma}

\begin{figure}[!ht]
\begin{center}
\includegraphics[width=0.25\linewidth]{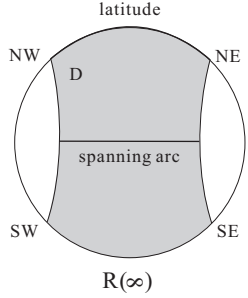}
\caption{A spanning arc and a latitude}
\label{spanningarc}
\end{center}
\end{figure}

Let $(B, t)$ be a tangle such that
$B \subset S^3 ( = \mathbb{R}^3\cup\{\infty\} )$ is
the complement of the unit $3$--ball $U$,
and $\partial t = \{ \mathrm{NW, NE, SE, SW} \}$.
We denote by $(B, t) + R(r)$ 
the knot or link in $S^3$ 
formed by the union of the strings of the tangles,
and let $\pi_r: X_r \to S^3 =B\cup U$ be
the two--fold cover branched along $(B, t) + R(r)$.

Suppose that $(B, t) + R(\infty)$ is a trivial knot. 
Then the two--fold cover $X_{\infty}$ branched along $(B, t) + R(\infty)$ is the $S^3$,
and the preimage of the spanning arc $\kappa$ for $R(\infty)$
is a knot in $X_{\infty} =S^3$,
which we call the \textit{covering knot} of $\kappa$. 
The exterior of the covering knot $K$ is $\pi_{\infty}^{-1}(B)$.
For $(B, t) +R(\infty)$
a replacement of $R(\infty)$ by a rational tangle $R(s)$
is called \textit{$s$--untangle surgery}
on $(B, t) +R(\infty)$. 
Performing untangle surgery downstairs
corresponds to replacing the solid torus $\pi_{\infty}^{-1}(U)$
by $\pi_s^{-1}(U)$ upstairs,
i.e.\ Dehn surgery on the covering knot $K$. 
This observation is referred to as the Montesinos trick \cite{Mon}. 
We denote the surgery slope by $\gamma_s$;
it is represented by a lift
of a meridian of $R(s)$. 
We say that $\gamma_s$ is the \textit{covering slope} of $s$. 
See the commutative diagram below. 

\begin{eqnarray*}
\begin{CD}
	S^3 @>\gamma_s\textrm{--surgery on } K >> K(\gamma) \\
@V{ \textrm{two--fold branched cover}}VV
		@VV{\textrm{two--fold branched cover}}V \\
	(B, t)\cup R(\infty) @>>s\textrm{--untangle surgery} > (B, t)\cup R(s)
\end{CD}
\end{eqnarray*}
\vskip 0.2cm
\begin{center}
\textsc{Diagram 2.} Montesinos trick 
\end{center}

For a link $L$ and an arc $\kappa$ with $\kappa \cap L = \partial \kappa$
we perform an untangle surgery along $\kappa$ as follows.
First take a regular neighborhood $N(\kappa)$ of $\kappa$ so that
$(N(\kappa), N(\kappa) \cap L)$ is a trivial tangle.
Then, identifying the trivial tangle $T =(N(\kappa), N(\kappa) \cap L)$
with the rational tangle $R(\infty)$,
we can replace $R(\infty) = T$ by a rational tangle $R(s)$;
this operation is called $s$--untangle surgery of $L$ along $\kappa$.
Note that the definition of $s$--untangle surgery along $\kappa$
relies on the identification of $T$ with $R(\infty)$.
If $L$ is a trivial knot, 
the two--fold cover of $S^3$ branched along $L$ is $S^3$,
and the preimage of $\kappa$ is a knot, which we call the covering knot
of $\kappa$.
Then, as before,
performing $s$--untangle surgery along $\kappa$ downstairs
corresponds to performing Dehn surgery on
the covering knot upstairs;
we call its surgery slope the covering slope. 

\bigskip

\section{Almost alternating unknots and L--space knots}
\label{construction}

A diagram of a knot is \textit{alternating} if over-crossings and under-crossings alternate while running along the diagram. 
A diagram of a knot is \textit{almost alternating} if the diagram is obtained by a single crossing change in an alternating diagram. 
Hence an almost alternating diagram has a crossing at which the crossing change makes the diagram alternating. 
Such a crossing is called a \textit{dealternator} of the almost alternating diagram. 
For later convenience, 
we call an arc connecting an over pass and an under pass at the dealternator a \textit{dealternating arc}. 
See Figure~\ref{dealternator}. 
(In the diagram, there are four dealternating arcs at each dealternator, 
but obviously they are isotopic.)

\begin{figure}[htbp]
\begin{center}
\includegraphics[width=0.28\linewidth]{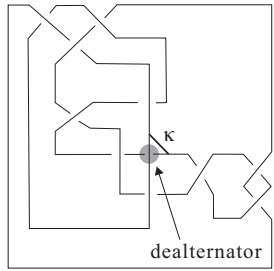}
\caption{$\kappa$ is a dealternating arc.}
\label{dealternator}
\end{center}
\end{figure}

The following result is useful to obtain an $L$--space knot, 
which was observed by Ozsv\'ath-Szab\'o \cite[8.3]{OS_unknotting}. 
For completeness we give its proof here.

\begin{theorem}[Ozsv\'ath-Szab\'o]
\label{OS}
Let $O$ be a trivial knot in $S^3$ and $\kappa$ an arc satisfying $\kappa \cap O = \partial \kappa$. 
Suppose that $O$ has an almost alternating diagram with $\kappa$ a dealternating arc.  
Then the covering knot $K$ of $\kappa$ is an L--space knot. 
\end{theorem}

\noindent
\textit{Proof of Theorem~\ref{OS}.}
First isotope $O \cup \kappa$ to a position so that the diagram of $O$ is almost alternating and $\kappa$ is a dealternating arc. 
Consider a tangle decomposition $(B, t) \cup R(\infty)$ of the trivial knot $O$ as depicted in Figure~\ref{crossing_untangle}(i). 
Then regard the dealternating arc $\kappa$ as a spanning arc and 
take the covering knot $K$ of $\kappa$. 
Figure~\ref{crossing_untangle} indicates that 
the crossing change at the dealternator (Figure~\ref{crossing_untangle}(ii)) 
corresponds to a $(-1/2)$--untangle surgery of $O$ along $\kappa$ (Figure~\ref{crossing_untangle}(iii)). 
Hence the $(-1/2)$--untangle surgery along $\kappa$ converts $O$ into an alternating knot $L$ whose alternating diagram is 
obtained from the almost alternating diagram of $O$ by crossing change at the dealternator. 
Since the two--fold branched cover of $S^3$ branched along a non-split alternating link is an L--space 
\cite[Lemma~3.2 and Proposition~3.3]{OSdoublecover}, 
$1/2$--surgery (in terms of $(\mu_{\infty}, \lambda)$--framing which may not be a preferred framing) 
on the covering knot $K$ yields an L--space. 
Thus $K$ is an L--space knot. 
\hspace*{\fill} $\square$Theorem~\ref{OS})

\medskip

\begin{figure}[htbp]
\begin{center}
\includegraphics[width=0.8\linewidth]{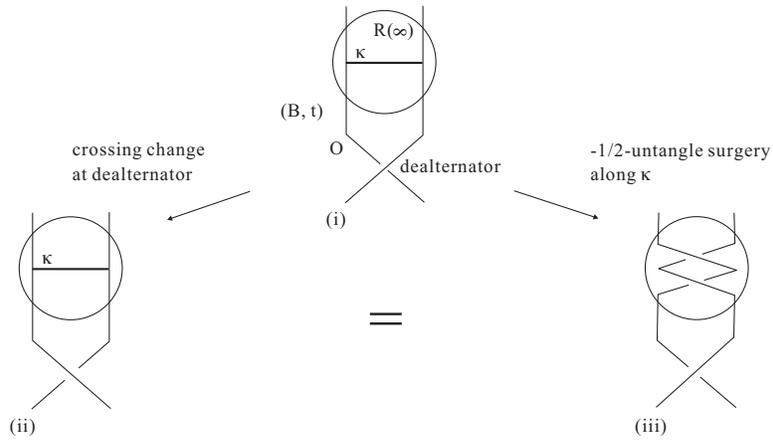}
\caption{Crossing change at the dealternator and $(-1/2)$--untangle surgery along the dealternating arc $\kappa$}
\label{crossing_untangle}
\end{center}
\end{figure}

\begin{remark}
\label{invertible}
By construction, 
L--space knots obtained in Theorem~\ref{OS} are strongly invertible. 
\end{remark}

Theorem~\ref{OS} shows that each almost alternating diagram of the unknot yields an L--space knot. 
So it is important to find such diagrams of the unknot. 
For this purpose, we recall a result of Tsukamoto \cite{Tsukamoto}. 

A \textit{reduced} diagram is one not containing any \textit{nugatory crossings} (Figure~\ref{nugatory}).  

\begin{figure}[htbp]
\begin{center}
\includegraphics[width=0.24\linewidth]{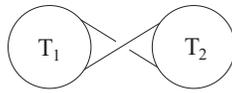}
\caption{A nugatory crossing}
\label{nugatory}
\end{center}
\end{figure}

Denote by $C_m$ a basic almost alternating diagram of the unknot depicted in Figure~\ref{Cm_moves}(i). 
Then Tsukamoto \cite{Tsukamoto} has shown that any reduced almost alternating diagram of the unknot can be obtained from $C_m$ 
by using only certain types of isotopies called \textit{flypes}, \textit{tongue moves} and \textit{twirl moves}; 
see Figure~\ref{Cm_moves}(ii). 
Recently McCoy \cite{Mc} has given an alternative proof of this result. 

\begin{theorem}[\cite{Tsukamoto, Mc}]
\label{almost alternating unknot} 
Any reduced almost alternating diagram of the unknot can be obtained from $C_m$ for some non-zero integer $m$, 
by a sequence of flypes, tongue moves and twirl moves. 
\end{theorem}

\begin{figure}[htbp]
\begin{center}
\includegraphics[width=0.9\linewidth]{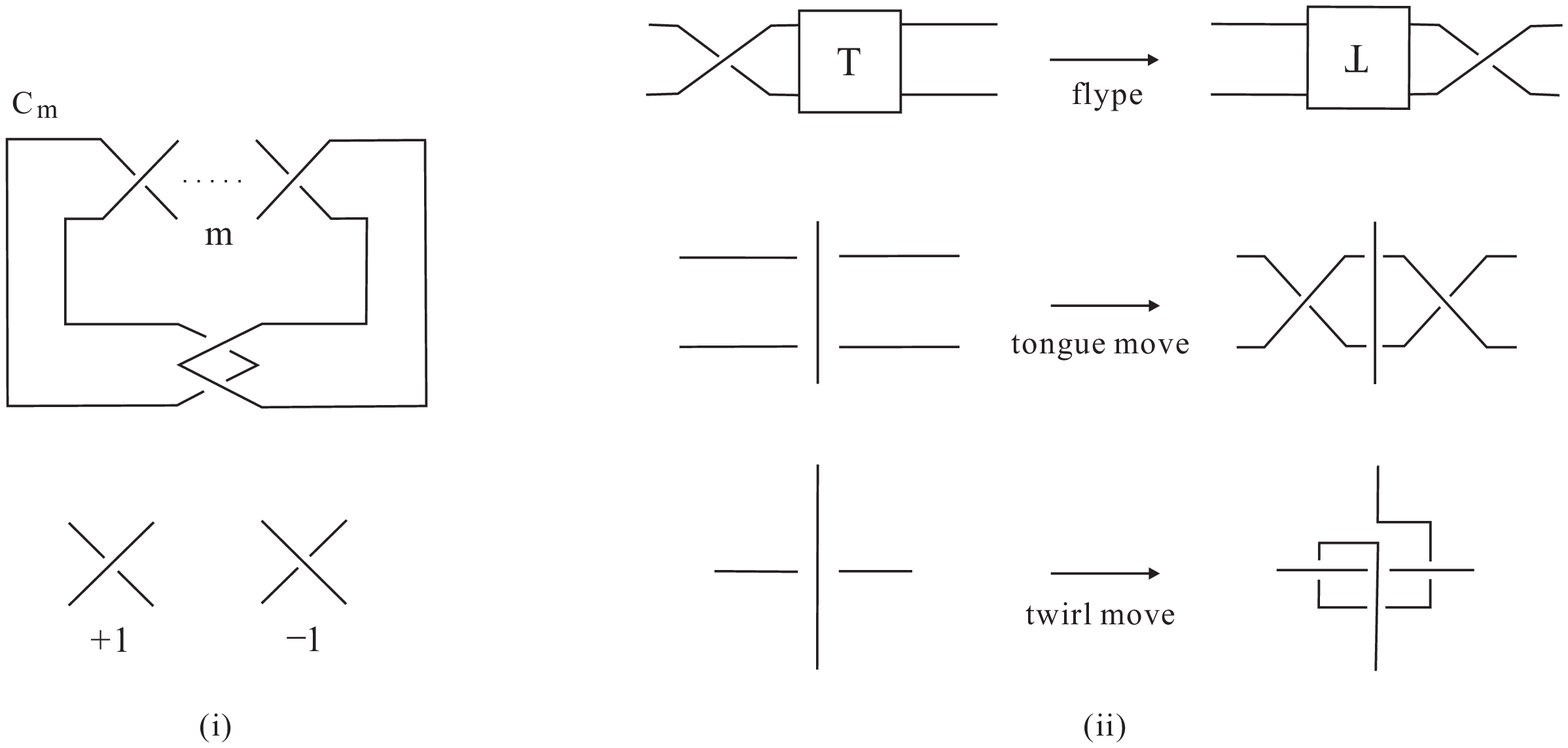}
\caption{A basic almost alternating diagram $C_m$ of the unknot, flype, tongue move and twirl move}
\label{Cm_moves}
\end{center}
\end{figure}

\section{A hyperbolic, L--space knot with no exceptional surgeries}

In this section we will prove Theorem~\ref{main} by giving a concrete example of 
a hyperbolic, L--space knot $K$ which has no exceptional surgeries, 
i.e. every nontrivial surgery on $K$ yields a hyperbolic $3$--manifold. 

\begin{theorem}
\label{L-space knot}
Let $K$ be a knot depicted in Figure~\ref{coveringknot}. 
Then $K$ is a hyperbolic, L--space knot which has no exceptional surgeries. 
In particular, it admits no Seifert fibered surgeries. 
\end{theorem}

\begin{figure}[htbp]
\begin{center}
\includegraphics[width=0.4\linewidth]{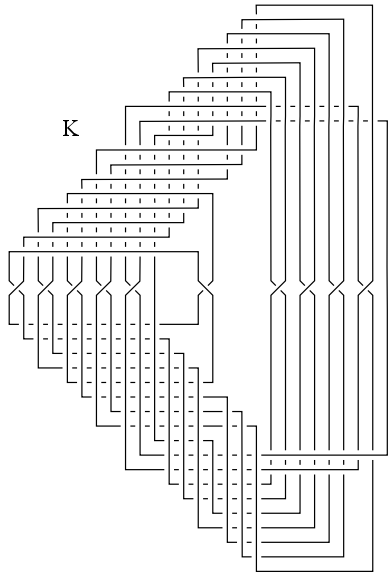}
\caption{A hyperbolic, L--space knot with no exceptional surgeries}
\label{coveringknot}
\end{center}
\end{figure}

\noindent
\textit{Proof of Theorem~\ref{coveringknot}.}
Based on Theorem~\ref{almost alternating unknot} we will apply flypes and tongue moves to 
$C_k$ several times to obtain a sufficiently complicated almost alternating diagram of the unknot. 
We start with the almost alternating diagram $C_{-3}$ of the unknot and apply a sequence of flypes and tongue moves as depicted in Figures~\ref{almostalt_unknot_1} and 
\ref{almostalt_unknot_2} to obtain the almost alternating diagram of the unknot given by the last picture of Figure~\ref{almostalt_unknot_2}. 

\begin{figure}[htbp]
\begin{center}
\includegraphics[width=1.0\linewidth]{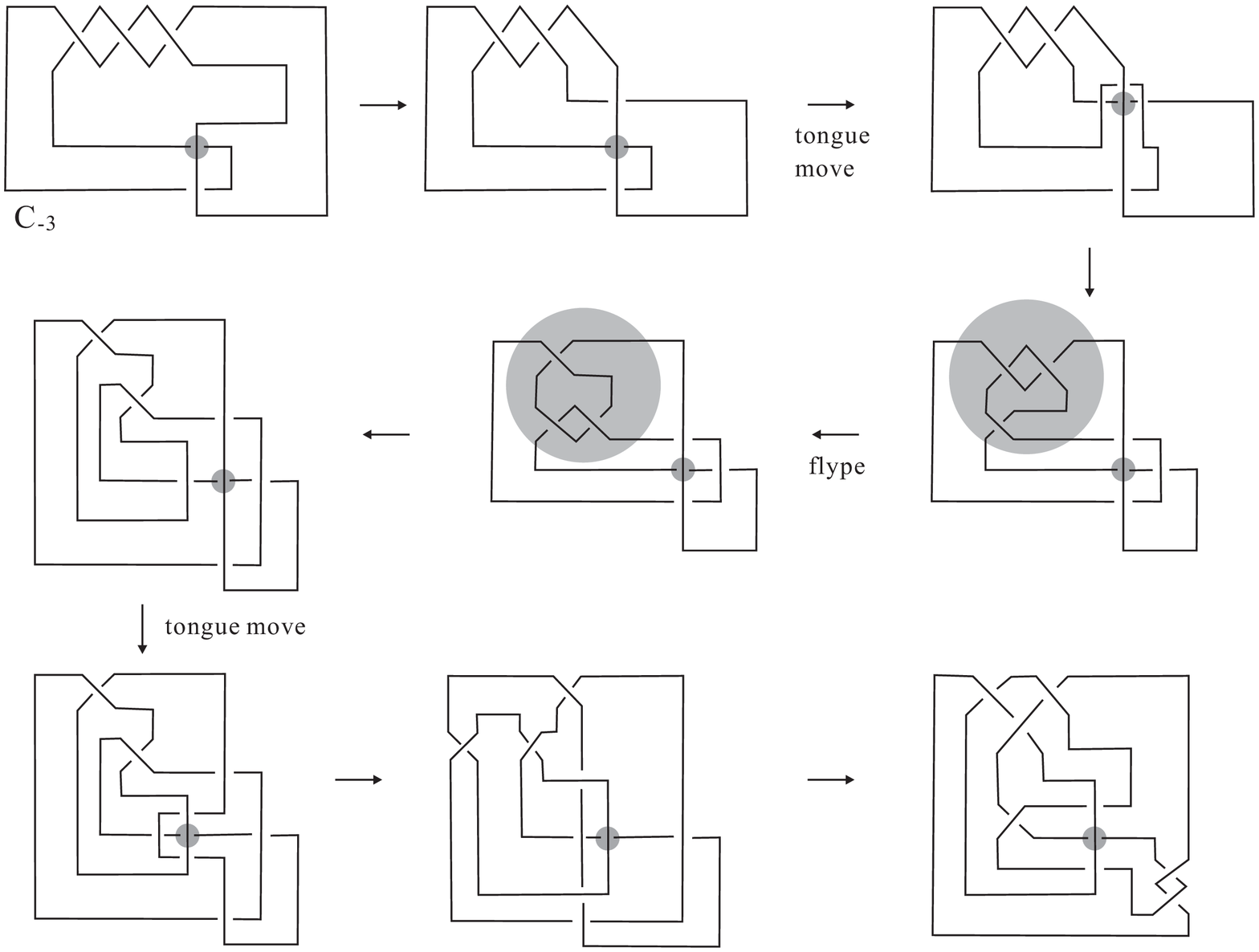}
\caption{Sequence of flypes and tongue moves}
\label{almostalt_unknot_1}
\end{center}
\end{figure}

\begin{figure}[htbp]
\begin{center}
\includegraphics[width=1.0\linewidth]{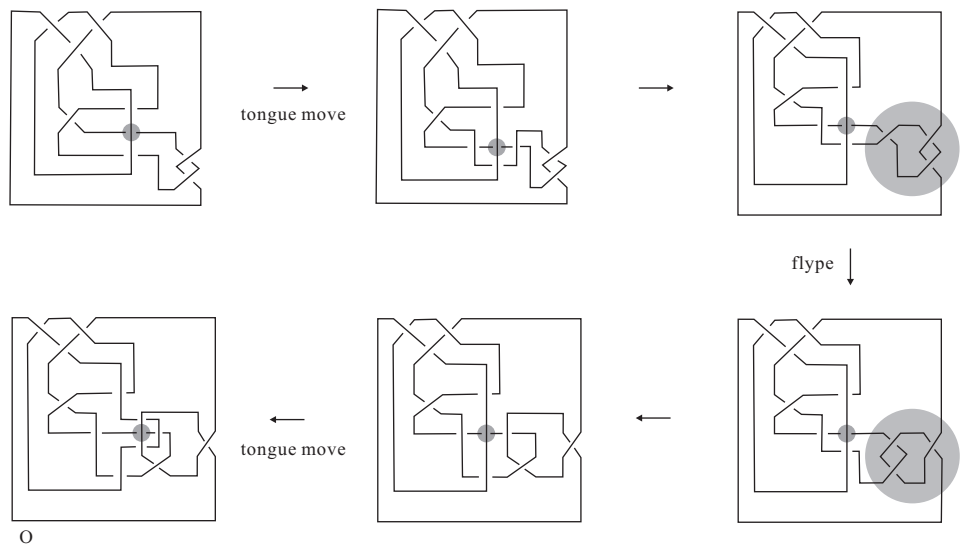}
\caption{Sequence of flypes and tongue moves; continued from Figure~\ref{almostalt_unknot_1}}
\label{almostalt_unknot_2}
\end{center}
\end{figure}

Let $\kappa$ be a dealternating arc as in Figure~\ref{coveringknot_isotopy1}. 
To obtain an explicit picture of the covering knot $K$ of $\kappa$, 
we apply isotopies given in Figures~\ref{coveringknot_isotopy1}--\ref{coveringknot_isotopy6}. 
Then taking the two--fold cover branched along $O$, 
we obtain the covering knot $K$ of $\kappa$; 
see Figure~\ref{coveringknot_branched}. 
By Theorem~\ref{OS} $K$ is an L--space knot.  

A proof showing that $K$ is a hyperbolic knot with no exceptional surgeries is 
computer-assisted but rigorous.  
SnapPy \cite{SnapPy} finds an approximated hyperbolic structure of $K$, 
which can be verified by the program HIKMOT \cite{HIKMOT}. 
To see that $K$ has no exceptional surgeries we run \texttt{fef.py} written by Ichihara and Masai \cite{IM}, 
which is a modification of a python code \texttt{find\_exceptional\_fillings.py} developed in \cite{MPR}. 
They improved the codes in \cite{MPR} using verified numerical analysis based on interval arithmetics to obtain 
mathematically rigorous computations. 
The code \texttt{fef.py}, together with HIKMOT, 
gives us as an output a list of candidates for exceptional fillings of $E(K)$, including all truly exceptional ones. 
See \cite[Section~6]{IM} for detailed explanation on \texttt{fef.py}. 
For the knot $K$, the set of candidate exceptional fillings turns out to be empty, 
and this proves that $K$ has no exceptional surgeries.  

\begin{figure}[htbp]
\begin{center}
\includegraphics[width=1.0\linewidth]{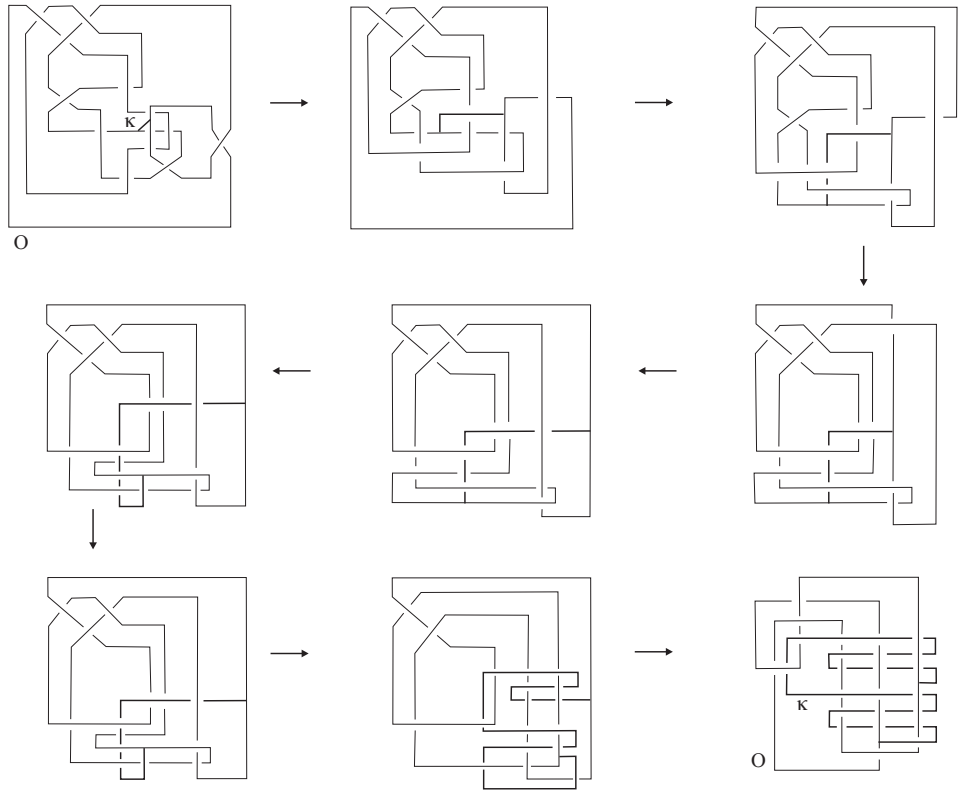}
\caption{Isotopies of $O \cup \kappa$}
\label{coveringknot_isotopy1}
\end{center}
\end{figure}

\begin{figure}[htbp]
\begin{center}
\includegraphics[width=1.0\linewidth]{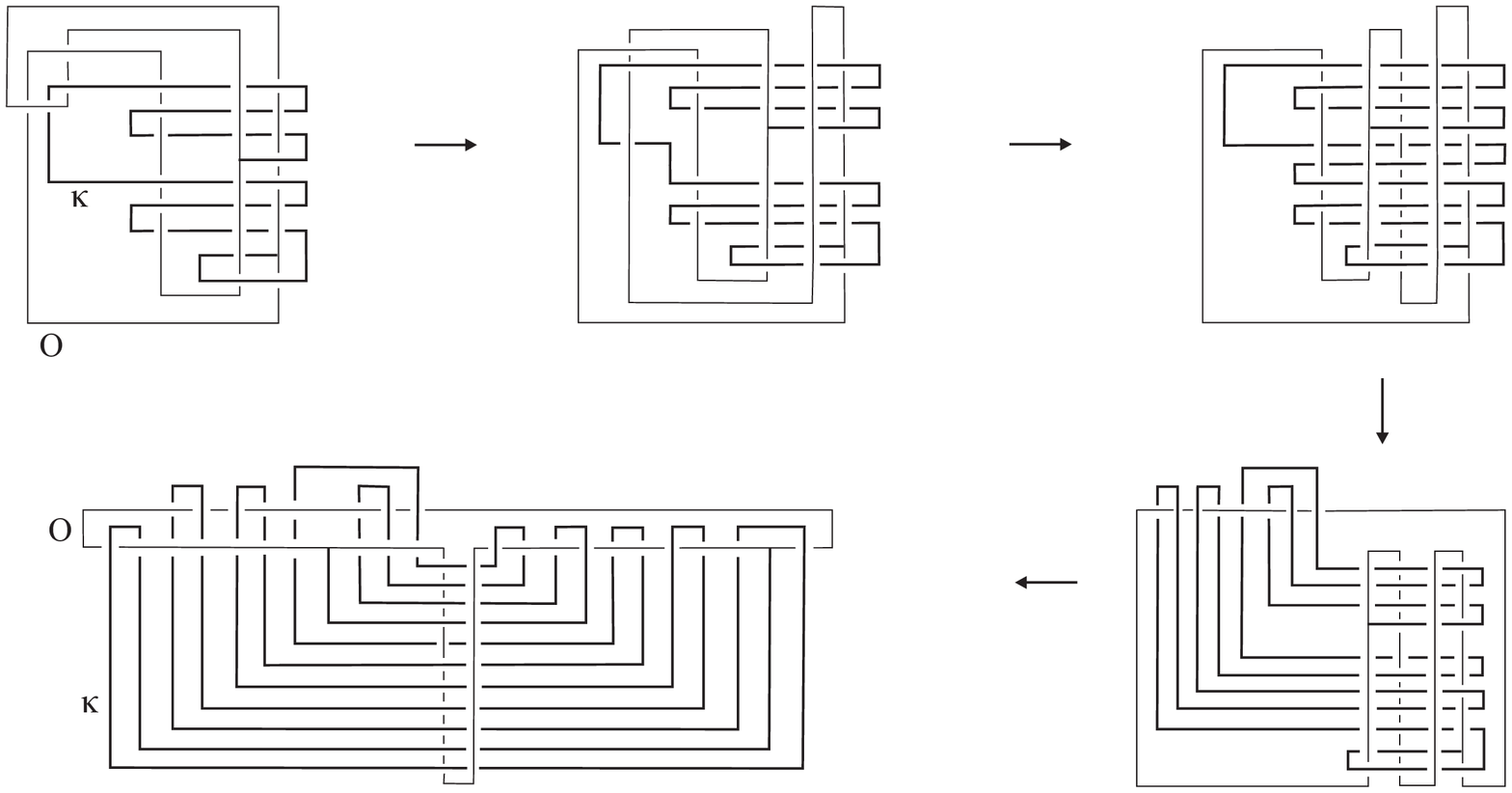}
\caption{Isotopies of $O \cup \kappa$; continued from Figure~\ref{coveringknot_isotopy1}}
\label{coveringknot_isotopy2}
\end{center}
\end{figure}

\begin{figure}[htbp]
\begin{center}
\includegraphics[width=1.0\linewidth]{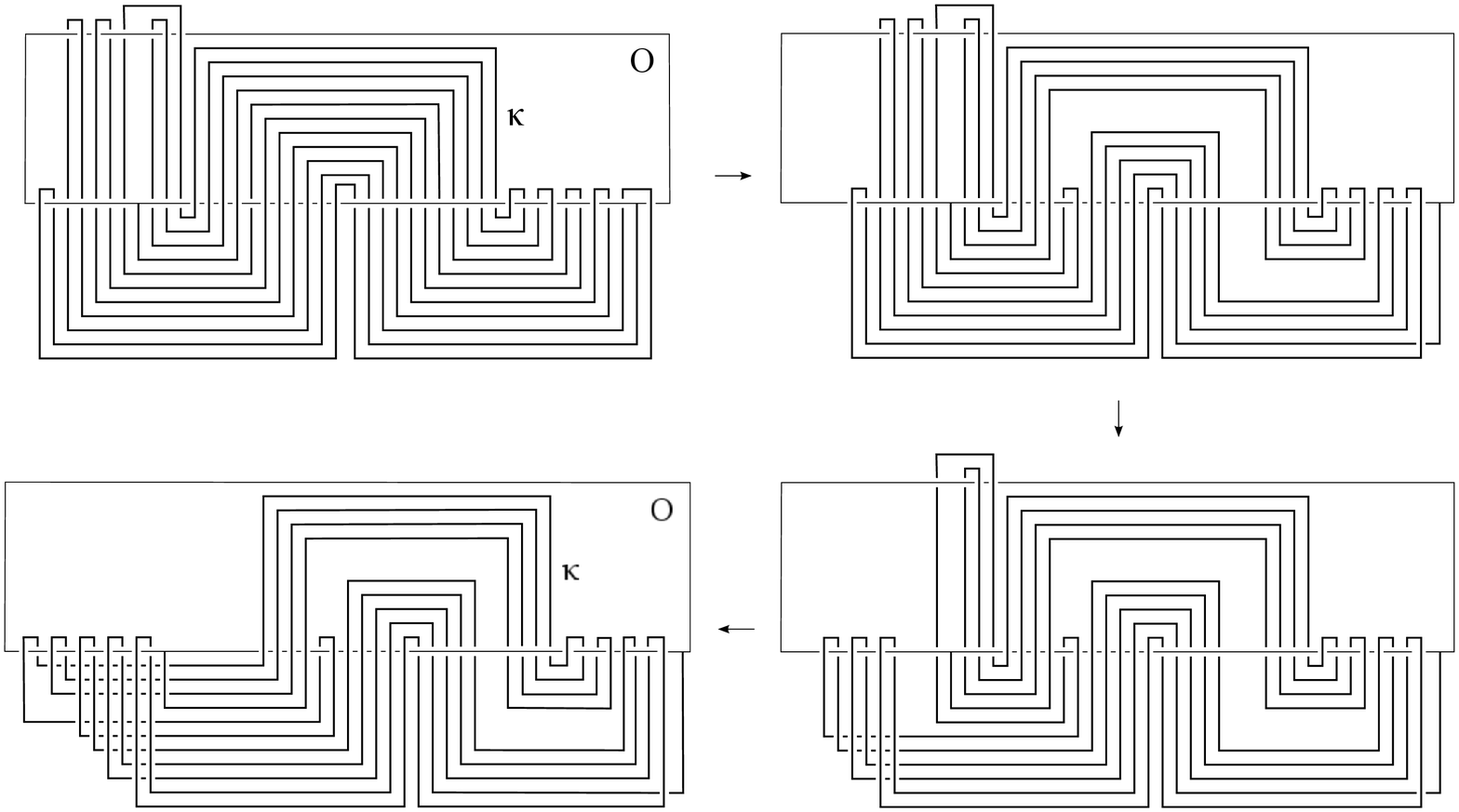}
\caption{Isotopies of $O \cup \kappa$; continued from Figure~\ref{coveringknot_isotopy2}}
\label{coveringknot_isotopy3}
\end{center}
\end{figure}

\begin{figure}[htbp]
\begin{center}
\includegraphics[width=1.0\linewidth]{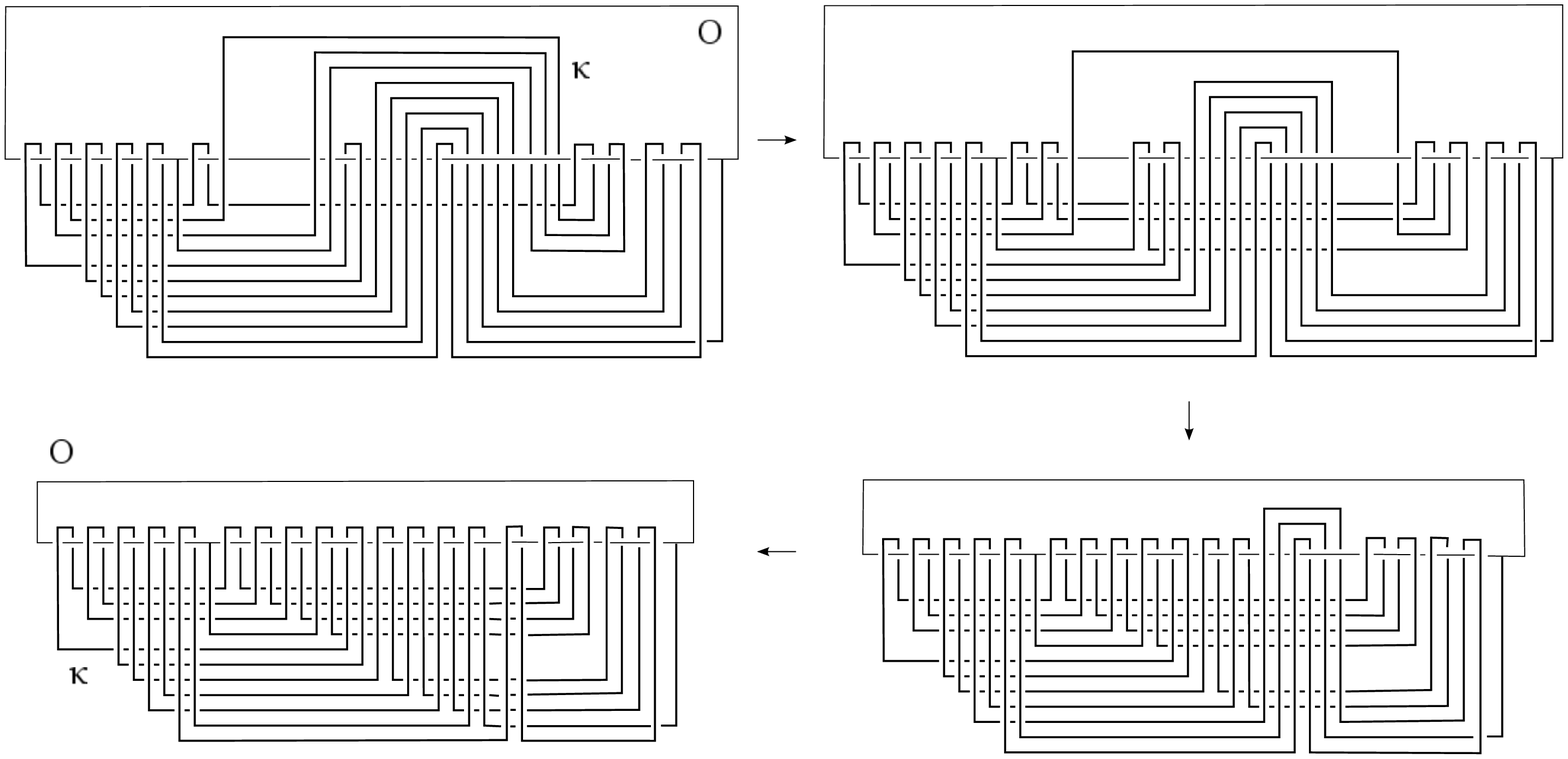}
\caption{Isotopies of $O \cup \kappa$; continued from Figure~\ref{coveringknot_isotopy3}}
\label{coveringknot_isotopy4}
\end{center}
\end{figure}

\begin{figure}[htbp]
\begin{center}
\includegraphics[width=1.0\linewidth]{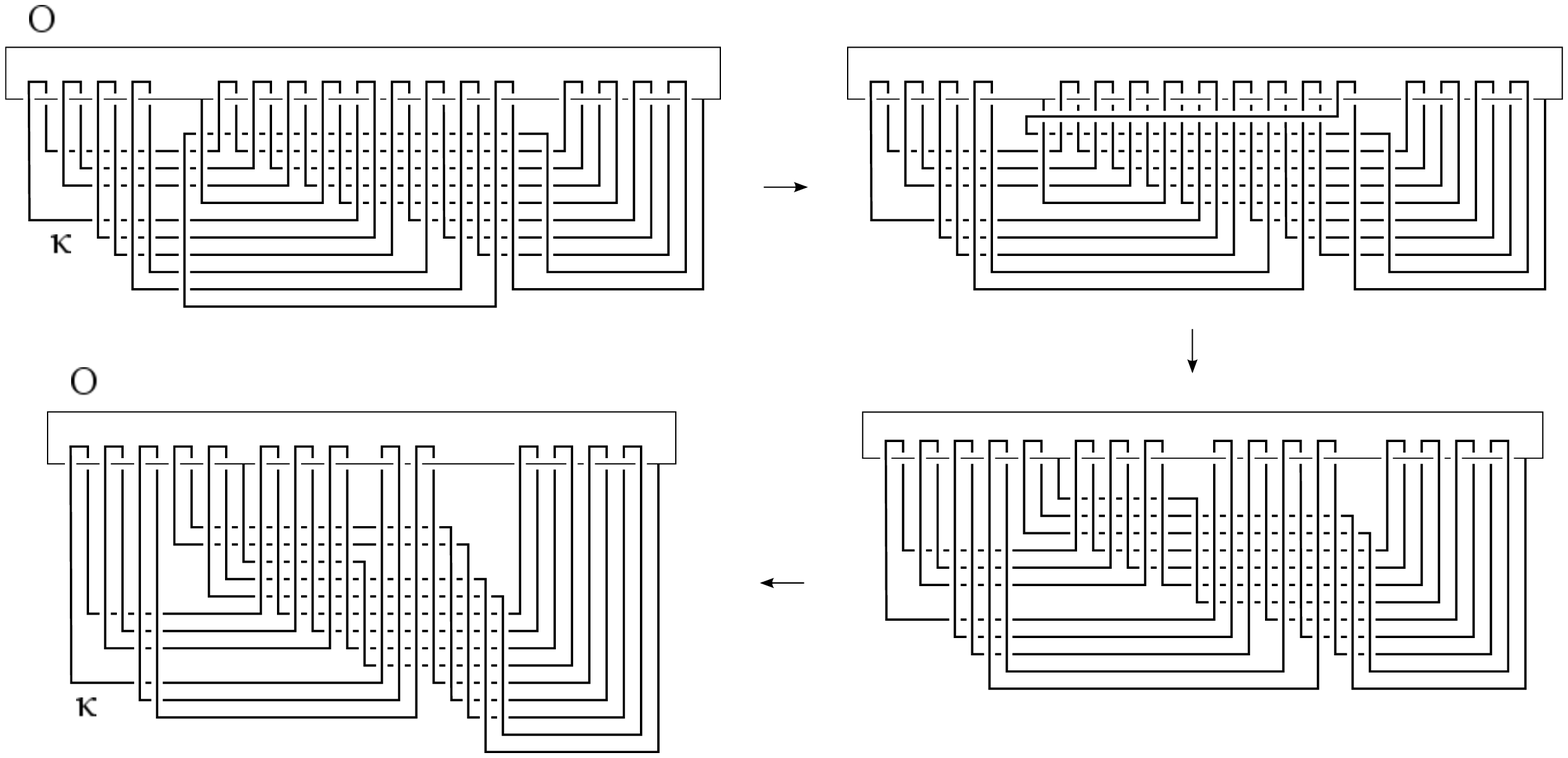}
\caption{Isotopies of $O \cup \kappa$; continued from Figure~\ref{coveringknot_isotopy4}}
\label{coveringknot_isotopy5}
\end{center}
\end{figure}

\begin{figure}[htbp]
\begin{center}
\includegraphics[width=1.0\linewidth]{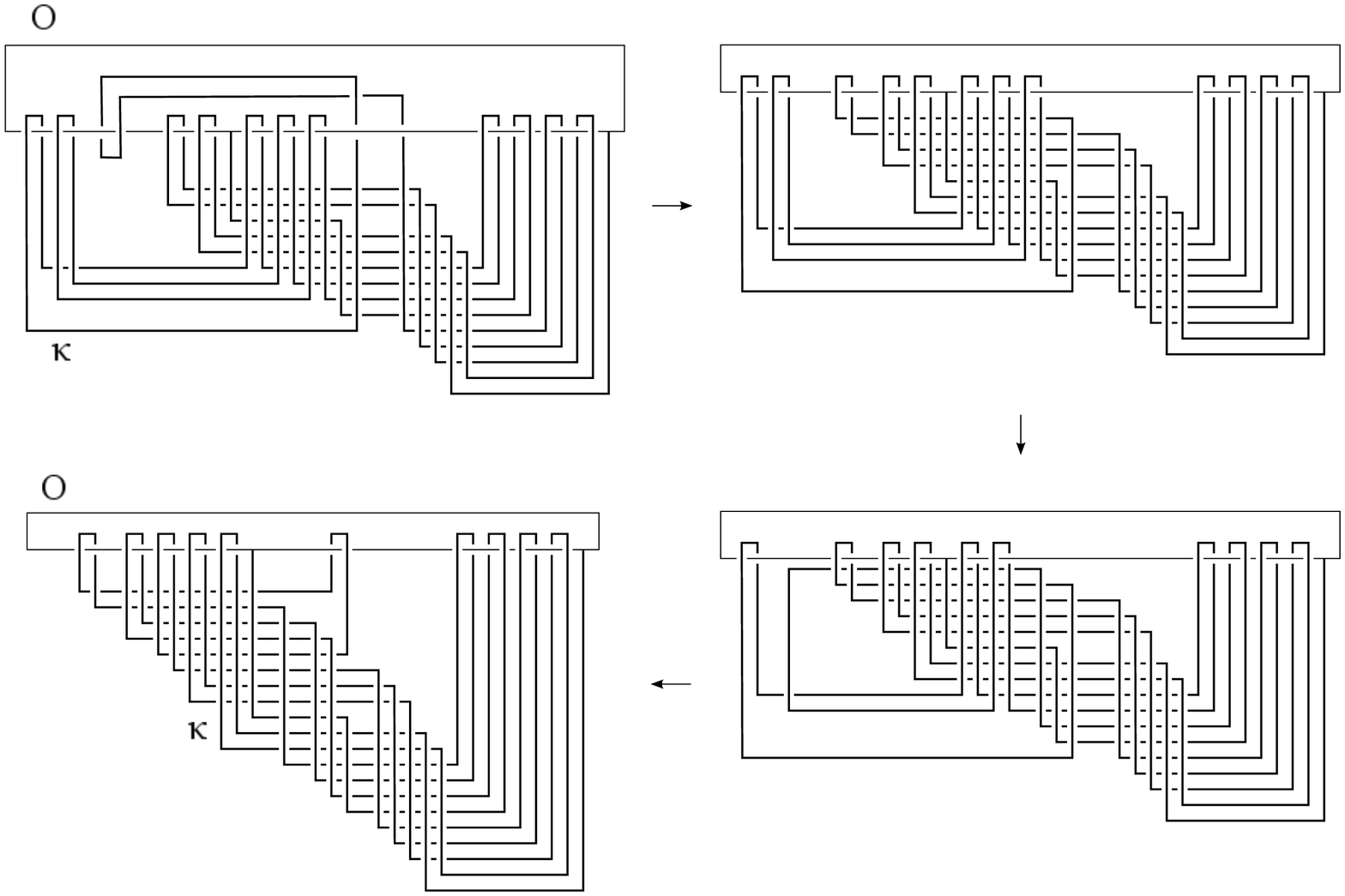}
\caption{Isotopies of $O \cup \kappa$; continued from Figure~\ref{coveringknot_isotopy5}}
\label{coveringknot_isotopy6}
\end{center}
\end{figure}

\begin{figure}[htbp]
\begin{center}
\includegraphics[width=1.0\linewidth]{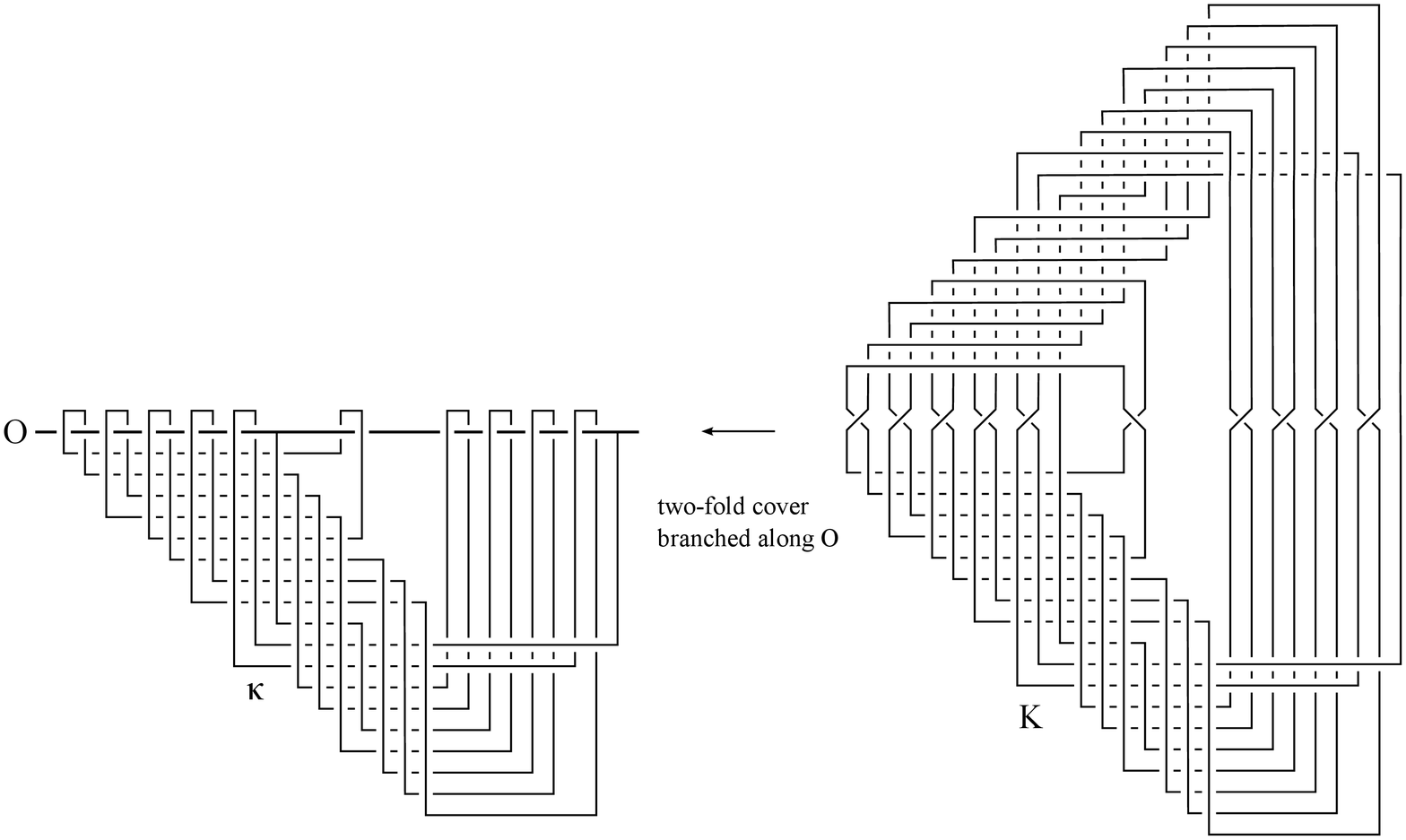}
\caption{The covering knot $K$ of $\kappa$}
\label{coveringknot_branched}
\end{center}
\end{figure}

\begin{remark}
\label{process}
In the fifth diagram in Figure~\ref{almostalt_unknot_1} 
if we perform $(-1)$--untangle surgery along the dealternating arc $\kappa_1$ given in Figure~\ref{Montesinos}(i), we
obtain a two-bridge knot. 
This implies that a surgery on the covering knot $K_1$ of $\kappa_1$ along the covering slope yields a lens space. 
Similarly, in the first diagram in Figure~\ref{almostalt_unknot_2}, 
performing $0$--untangle surgery along the dealternating arc $\kappa_2$ given in Figure~\ref{Montesinos}(ii) yields also a two-bridge knot. 
Thus the covering knot $K_2$ of $\kappa_2$ has a lens space surgery.  
\end{remark}

\begin{figure}[htbp]
\begin{center}
\includegraphics[width=0.7\linewidth]{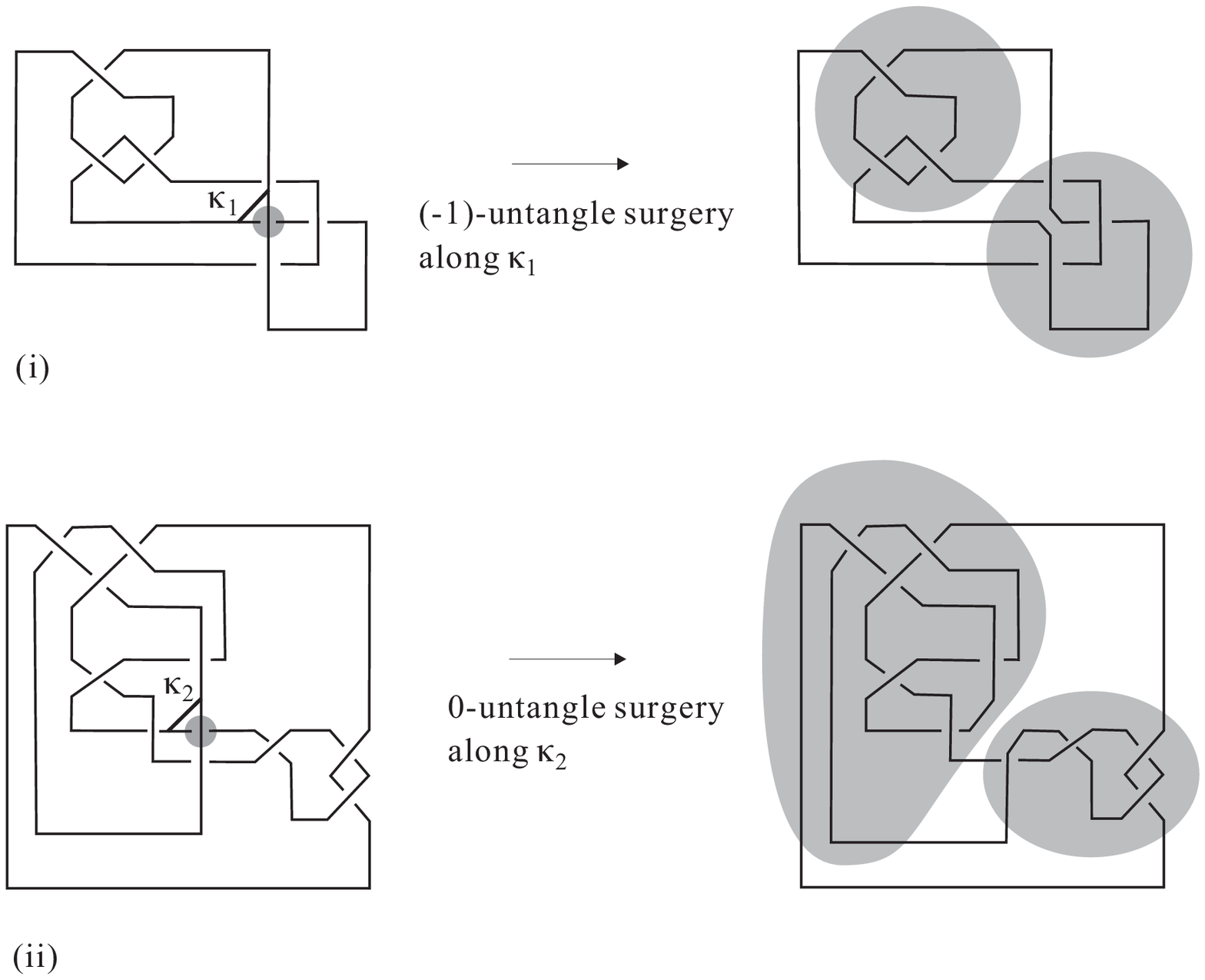}
\caption{$(-1)$--untangle surgery along $\kappa_1$ yields a two--bridge knot, 
and $0$--untangle surgery along $\kappa_2$ yields a two--bridge knot. }
\label{Montesinos}
\end{center}
\end{figure}


\bigskip


\begin{thebibliography}{99}

\bibitem{BM}
K. Baker and A. Moore; 
Montesinos knots, Hopf plumbings, and L--space surgeries, 
preprint 2014. 

\bibitem{BePe}
R. Benedetti and C. Petronio; 
Lectures on hyperbolic geometry, 
Universitext, Springer-Verlag, 1992. 

\bibitem{Berge2} 
J. Berge; 
Some knots with surgeries yielding lens spaces, 
unpublished manuscript.

\bibitem{BoileauPorti}
M. Boileau and J. Porti; 
Geometrization of 3-orbifolds of cyclic type,  
Ast\'erisque \textbf{272} (2001), 208pp. 

\bibitem{SnapPy}
M. Culler, N. Dunfield, and J. R. Weeks; 
SnapPy, a computer program for studying the
geometry and topology of 3-manifolds, \texttt{http://snappy.computop.org}. 

\bibitem{Hedden}
M. Hedden; 
On knot Floer homology and cabling II,  
Int.\ Math.\ Res.\ Not.\ IMRN, (12):2248--2274,
2009.

\bibitem{HIKMOT}
N. Hoffman, K. Ichihara, M. Kashiwagi, H. Masai, S. Oishi, and A. Takayasu; 
Verified computations for hyperbolic 3-manifolds, preprint 2013. 

\bibitem{HLV}
J. Hom, T. Lidman and F. Vafaee; 
Berge-Gabai knots and L--space satellite operations, preprint 2014. 

\bibitem{IM}
K. Ichihara and H. Masai; 
Exceptional surgeries on alternating knots, preprint 2014. 

\bibitem{LM}
T. Lidman and A. Moore; 
Pretzel knots with L--space surgeries, 
preprint 2013. 

\bibitem{MPR}
B. Martelli, C. Petronio, and F. Roukema; 
Exceptional Dehn surgery on the minimally twisted five-chain link, 
Comm.\ Anal.\ Geom.\ \textbf{22} (2014) 689--735. 

\bibitem{Mc}
D. McCoy; 
Alternating knots with unknotting number one, 
preprint 2014. 

\bibitem{Mon}
J. M. Montesinos; 
Surgery on links and double branched coverings of $S^3$, 
Knots, groups, and $3$--manifolds (Papers dedicated to the memory of R.H.Fox),
Ann. Math. Studies, 84, Princeton Univ. Press 1975, 227--259,  

\bibitem{Mote}
K. Motegi; 
L--space surgery and twisting operation, 
preprint 2014. 

\bibitem{OS_unknotting}
P. Ozsv\'ath and Z. Szab\'o; 
Knots with unknotting number one and Heegaard Floer homology, 
Topology \textbf{44} (2005), 705--745. 

\bibitem{OS3} 
P. Ozsv\'ath and Z. Szab\'o; 
On knot Floer homology and lens space surgeries, 
Topology \textbf{44} (2005), 1281--1300. 

\bibitem{OSdoublecover}
P. Ozsv\'ath and Z. Szab\'o; 
On the Heegaard Floer homology of branched double-covers. 
Adv.\ Math.\ 
\textbf{194} (2005), 1--33. 

\bibitem{OS4} 
P. Ozsv\'ath and Z. Szab\'o; 
Knot Floer homology and rational surgeries, 
Algebr.\ Geom.\ Topol.\ \textbf{11} (2011), 1--68. 


\bibitem{PetPorti}
C. Petronio and J. Porti;
Negatively oriented ideal triangulations and a proof of Thurston's hyperbolic Dehn filling theorem, 
Expo.\ Math.\ \textbf{18} (2000), 1--35. 


\bibitem{T1}
W. P. Thurston; 
The geometry and topology of $3$-manifolds, 
Lecture notes, Princeton University, 1979. 

\bibitem{T2}
W. P. Thurston; 
Three dimensional manifolds, Kleinian groups and hyperbolic geometry, 
Bull.\ Amer.\ Math.\ Soc.\ \textbf{6} (1982), 357--381. 

\bibitem{Tsukamoto}
T. Tsukamoto; 
The almost alternating diagram of the trivial knot, 
J.\ Topology \textbf{2} (2009), 77--104. 


\end{thebibliography}
\end{document}